\documentclass[11pt]{article}
\usepackage{enumerate}
\usepackage{amssymb,a4wide,latexsym,makeidx,epsfig,fleqn}
\usepackage{amsthm}
\usepackage{amsmath}
\usepackage{enumerate}
\newtheorem{theorem}{Theorem}[section]

\newtheorem{lemma}[theorem]{Lemma}

\newtheorem{corollary}[theorem]{Corollary}

\begin{document}
\textwidth 150mm \textheight 225mm
\title{The signless Laplacian spectral radius of  subgraphs of regular graphs
\thanks{ Supported by the National Natural Science Foundation of China (No. 11171273)
 \vskip 0.05in}}
\author{{Qi Kong$^1$, Ligong Wang$^2$}\\
{\small Department of Applied Mathematics, School of Science, Northwestern
Polytechnical University,}\\ {\small  Xi'an, Shaanxi 710072,
People's Republic
of China.}\\ {\small $^1$E-mail: kongqixgd@163.com}\\
{\small $^2$E-mail: lgwangmath@163.com}}
\maketitle
\begin{center}
\begin{minipage}{120mm}
\vskip 0.3cm
\begin{center}
{\small {\bf Abstract}}
\end{center}
{\small Let $q(H)$ be the signless Laplacian spectral radius of a graph $H$. In this paper, we prove that
\\1. Let $H$ be a proper subgraph of a $\Delta$-regular graph $G$ with $n$ vertices and diameter $D$. Then
$$2\Delta - q(H)>\frac{1}{n(D-\frac{1}{4})}.$$
\\2. Let $H$ be a proper subgraph of a $k$-connected $\Delta$-regular graph $G$ with $n$ vertices, where $k\geq 2$. Then
$$2\Delta-q(H)>\frac{2(k-1)^{2}}{2(n-\Delta)(n-\Delta+2k-4)+(n+1)(k-1)^{2}}.$$
Finally, we compare the two bounds. We obtain that when $k>2\sqrt{\frac{(n-\Delta)(n+\Delta-4)}{n(4D-3)-2}}+1$, the second bound is always better than the first. On the other hand, when $k<\frac{2(n-\Delta)}{\sqrt{n(4D-3)-2}}+1$, the first bound is always better than the second.

 }
\end{minipage}
\end{center}

\section{Introduction }
\label{sec:ch6-introduction}

As usual, let $G=(V(G), E(G))$ be a finite, undirected and simple graph with vertex set $V(G)=\{v_{1}, v_{2}, \cdots, v_{n}\}$ and edge set $E(G)$. Set $N_{G}(v_{i})=\{v|v_{i}v\in E(G)\}$ and $d_{G}(v_{i})=|N_{G}(v_{i})|$, or simply $N(v_{i})$ and $d_{i}=d(v_{i})$, respectively. Let $\delta=\delta(G)$ and $\Delta=\Delta(G)$ denote the minimum degree and maximum degree of the graph $G$, respectively. If $\Delta=\delta$, then $G$ is regular. Let $P=x_{1}x_{2}\ldots x_{t}$ be a path in $G$ with a given orientation. We denote by $x_{i}Px_{j}$ the path $x_{i}x_{i+1}\ldots x_{j-1}x_{j}$ for $i<j$. The distance between any two vertices $v_{i}$ and $v_{j}$ in $G$ is the number of edges in a shortest path connecting $v_{i}$ and $v_{j}$, denoted by $d_{G}(v_{i}, v_{j})$. The diameter $D=D_{G}$ of $G$ is the maximum distance between any two vertices of $G$. The (vertex) connectivity $\kappa(G)$ of $G$ is the minimum number of vertices whose removal disconnects $G$ or reduces it to a single vertex. For an integer $k\geq 1$, $G$ is called $k$-connected if $\kappa(G)\geq k$. For terminologies and notations of graphs undefined here, we refer the reader to \cite{MGT}.

Let $A(G)$ be the adjacency matrix of $G$ and $D(G)=diag(d_{1}, d_{2}, \cdots, d_{n})$ be the diagonal matrix of vertex degrees of $G$. The matrix $L(G)=D(G)-A(G)$ is called the Laplacian matrix of $G$, and the matrix $Q(G)=D(G)+A(G)$ is called the signless Laplacian matrix of $G$. The largest eigenvalue of $A(G)$, $L(G)$ and $Q(G)$ are called spectral radius, Laplacian spectral radius and signless Laplacian spectral radius of $G$, and denoted by $\rho(G)$, $\mu(G)$ and $q(G)$, respectively. Since $A(G)$, $L(G)$ and $Q(G)$ are real symmetric matrices, their eigenvalues are real numbers.

If $G$ is a simple connected graph, then the matrix $A(G)$ (or $Q(G)$) is a nonnegative irreducible matrix and the largest eigenvalues of $A(G)$ (or $Q(G)$) is nonnegative. By Perron-Frobenius Theorem, $\rho(G)$ (or $q(G)$) is simple and has a unique positive unit eigenvector.

We all know that $\rho(G)\leq \Delta(G)$ with equality if and only if $G$ is regular. Some good bounds on the spectral radius $\rho(G)$ of connected irregular graphs have been obtained by various authors in \cite{TSR, EEO, TEE, OTL, TLE, LTSR, EAE}.
Moreover, if $H$ is a proper subgraph of a connected graph $G$, then $\rho(G)>\rho(H)$.  Then Nikiforov in \cite{RTC} gave a bound of $\rho(G)-\rho(H)$.
So combining the above two famous results, the authors in \cite{TEE, VTSR, RTC, LTSR} obtained some bounds of $\rho(G)-\rho(H)$ when $H$ is the proper subgraph of connected regular graph $G$.

Also, as we all know that $q(G)\leq 2\Delta(G)$ with equality if and only if $G$ is regular \cite{AIT}. In fact, in 2013, Ning et al. \cite{OTS} gave a bound on the signless Laplacian spectral radius of irregular graph $G$ with $n$ vertices, maximum degree $\Delta$ and diameter $D$:
\begin{equation}\label{eq:e1}
2\Delta-q(G)>\frac{1}{n(D-\frac{1}{4})}.
\end{equation}
And in 2015, Chen and Hou \cite{OTB} obtained a bound on the signless Laplacian spectral radius of $k$-connected irregular graph $G$:
\begin{equation}\label{eq:e2}
2\Delta - q(G)>\frac{2(n\Delta-2m)k^{2}}{2(n\Delta-2m)[n^{2}-2(n-k)]+nk^{2}}.
\end{equation}
In \cite{OTB}, they also obtained when $k\geq \sqrt{n}$, bound (\ref{eq:e2}) is always better than bound (\ref{eq:e1}).

We also know that $q(H)\leq q(G)$ whenever $H$ is a subgraph of $G$. So we can arise the following question:
\begin{center}How small $q(G)-q(H)$ can be when $H$ is a subgraph of a regular graph $G$?
\end{center}

In this paper, we give two bounds of $q(G)-q(H)$ when $H$ is a subgraph of a regular graph $G$.
\noindent\begin{theorem}\label{th:t1}
Let $H$ be a proper subgraph of a $\Delta$-regular graph $G$ with $n$ vertices and diameter $D$. Then
\begin{equation}\label{eq:e3}
2\Delta - q(H)>\frac{1}{n(D-\frac{1}{4})}.
\end{equation}
\end{theorem}
By taking connectivity parameter into account, we establish the following theorem.
\noindent\begin{theorem}\label{th:t2}
Let $H$ be a proper subgraph of a $k$-connected $\Delta$-regular graph $G$ with $n$ vertices. If $k\geq 2$, then
\begin{equation}\label{eq:e4}
2\Delta-q(H)> \frac{2(k-1)^{2}}{2(n-\Delta)(n-\Delta+2k-4)+(n+1)(k-1)^{2}}.
\end{equation}
\end{theorem}
Finally, we compare the two bounds. We also obtain when $k>2\sqrt{\frac{(n-\Delta)(n+\Delta-4)}{n(4D-3)-2}}+1$, bound (\ref{eq:e4}) is always better than bound (\ref{eq:e3}). On the other hand, when $k<\frac{2(n-\Delta)}{\sqrt{n(4D-3)-2}}+1$, bound (\ref{eq:e3}) is always better than bound (\ref{eq:e4}).

Moreover, we notice that $\mu(G)\leq q(G)$ when $G$ is a graph, and if $G$ is connected, then the equality holds if and only if $G$ is a bipartite graph \cite{ZTSR}. Then we can give two upper bound of Laplacian spectral radius of subgraphs of regular graphs.
\noindent\begin{corollary}\label{co:c1}
Let $H$ be a proper subgraph of a $\Delta$-regular graph $G$ with $n$ vertices and diameter $D$. Then
$$2\Delta - \mu(H)>\frac{1}{n(D-\frac{1}{4})}.$$
\end{corollary}
\noindent\begin{corollary}\label{co:c2}
Let $H$ be a proper subgraph of a $k$-connected $\Delta$-regular graph $G$ with $n$ vertices. If $k\geq 2$, then
$$2\Delta-\mu(H)> \frac{2(k-1)^{2}}{2(n-\Delta)(n-\Delta+2k-4)+(n+1)(k-1)^{2}}.$$
\end{corollary}

\section{The proofs of Theorems \ref{th:t1}  and \ref{th:t2}}
In this section, we begin to prove Theorems \ref{th:t1} and \ref{th:t2}. Before our proofs we give a lemma which is used in the proofs. It is an immediate consequence of the Cauchy-Schwarz inequality (or see \cite{LTSR}).
\noindent\begin{lemma}\label{le:l2}(\cite{LTSR})
If $a, b>0$, then $a(x-y)^{2}+by^{2}\geq abx^{2}/(a+b)$ with equality if and only if $y=ax/(a+b)$.
\end{lemma}
\textit{Proof of Theorem \ref{th:t1}.}  Let $G$ be a $\Delta$-regular graph. And we suppose that $H$ is a maximal proper subgraph of $G$, i.e., $V(H)=V(G)$ and $H$ differs from $G$ in a single edge $uv$, i.e., $H=G-uv$. Then $d_{H}(u)=d_{H}(v)=\Delta-1$.

Let $x=(x_{1}, x_{2}, \cdots, x_{n})^{T}$ be the unique unit positive eigenvector of $Q(H)$ corresponding to $q(H)$. Clearly, $x_{1}^{2}+x_{2}^{2}+\cdots+x_{n}^{2}=1$. Let $w$ be a vertex such that $x_{w}=max_{1\leq i\leq n}x_{i}$. Thus we have $x_{w}> \frac{1}{\sqrt{n}}$.

We will prove that $u\neq w$ and $v\neq w$. Indeed, if $u=w$, then
$$q(H)x_{u}=(\Delta-1)x_{u}+\sum\limits_{uv_{i}\in E(H)}x_{i}\leq 2(\Delta-1)x_{u},$$
and thus $q(H)\leq 2(\Delta-1)$, contradicting the fact that $q(H)> 2\delta=2(\Delta-1)$ since $H$ is a irregular graph (see \cite{AIT}). Hence, $u\neq w$. Similarly, $v\neq w$.

We also find that
\begin{align*}
\displaystyle 2\Delta-q(H)&=2\Delta\cdot 1-x^{T}Q(H)x\\
&=2\Delta\sum\limits_{i=1}^{n}x_{i}^{2}-\sum\limits_{i=1}^{n}d_{i}x_{i}^{2}-2\sum\limits_{v_{i}v_{j}\in E(H)}x_{i}x_{j}\\
            &=2(x_{u}^{2}+x_{v}^{2})+\sum\limits_{i=1}^{n}d_{i}x_{i}^{2}-2\sum\limits_{v_{i}v_{j}\in E(H)}x_{i}x_{j}\\
            &=2(x_{u}^{2}+x_{v}^{2})+\sum\limits_{v_{i}v_{j}\in E(H)}(x_{i}-x_{j})^{2}.
 \end{align*}
Next we consider the following two cases.

{\it Case 1.} $d_{H}(w, u)\leq D-1$.

Select a shortest path $u=u_{0},u_{1},\ldots ,u_{l}=w$ joining $u$ to $w$ in $H$, i.e., $l\leq D-1$. By Lemma \ref{le:l2} and the Cauchy-Schwarz inequality, we have
\begin{align*}
\displaystyle 2\Delta-q(H)&=2(x_{u}^{2}+x_{v}^{2})+\sum\limits_{v_{i}v_{j}\in E(H)}(x_{i}-x_{j})^{2}\\
&> \sum\limits_{i=0}^{l-1}(x_{u_{i}}-x_{u_{i+1}})^{2}+2x_{u}^{2}
\end{align*}
\begin{align*}&\geq \frac{1}{l}[\sum\limits_{i=0}^{l-1}(x_{u_{i}}-x_{u_{i+1}})]^{2}+2x_{u}^{2}\\
&=\frac{1}{l}(x_{w}-x_{u})^{2}+2x_{u}^{2}\\
&\geq \frac{2}{2l+1}x_{w}^{2}>\frac{2}{2(D-1)+1}\cdot\frac{1}{n}\\
&=\frac{1}{n(D-\frac{1}{2})}>\frac{1}{n(D-\frac{1}{4})}.
 \end{align*}
{\it Case 2.} $d_{H}(w, u)\geq D$.

 In this case, by symmetry, $d_{H}(w, v)\geq D$. Let $P: u=u_{0}, u_{1}, \ldots, u_{l}=w$ and $Q$ be shortest paths joining $u$ to $w$ and $v$ to $w$ in $G$, respectively. Next we will prove that $u\notin Q$ and $v\notin P$.

If $u\in Q$, then there exists a path of length at most $D-1$ joining $w$ to $u$ in $G$, and thus in $H$, a contradiction. Hence, $u\notin Q$. By symmetry, $v\notin P$.

Thus the paths $P$ and $Q$ belong to $H$, and we have
$$d_{H}(w, u)=d_{H}(w, v)=D.$$
Then we have $l=D$. Let $t$ be the smallest index $j$ such that $u_{j}$ is on $Q$, then $t\geq 1$. Obviously, $uPu_{t}$ and $vQu_{t}$ have the same length. Using Lemma \ref{le:l2} and the Cauchy-Schwarz inequality, it follows that
\begin{align*}
\displaystyle 2\Delta-q(H)&= 2(x_{u}^{2}+x_{v}^{2})+\sum\limits_{v_{i}v_{j}\in E(H)}(x_{i}-x_{j})^{2}\\
&\geq 2(x_{u}^{2}+x_{v}^{2})+\sum\limits_{i=0}^{t-1}(x_{u_{i}}-x_{u_{i+1}})^{2}+\sum\limits_{ij\in E(vQu_{t})}(x_{i}-x_{j})^{2}\\
&~~~~~+\sum\limits_{i=t}^{D-1}(x_{u_{i}}-x_{u_{i+1}})^{2}\\
&\geq \frac{1}{t}(x_{u_{t}}-x_{u})^{2}+2x_{u}^{2}+\frac{1}{t}(x_{u_{t}}-x_{v})^{2}+2x_{v}^{2}
+\frac{1}{D-t}(x_{w}-x_{u_{t}})^{2}\\
&\geq \frac{4}{2t+1}x_{u_{t}}^{2}+\frac{1}{D-t}(x_{w}-x_{u_{t}})^{2}\\
&\geq\frac{4}{4D-2t+1}x_{w}^{2}\geq\frac{4}{4D-1}x_{w}^{2}\\
&> \frac{1}{n(D-\frac{1}{4})}.
 \end{align*}
This completes the  proof of Theorem \ref{th:t1}. \quad $\square$

\textit{Proof of Theorem \ref{th:t2}.}
Let $G$ be a $k$-connected $\Delta$-regular graph. And we suppose that $H$ is a maximal proper subgraph of $G$, i.e., $V(H)=V(G)$ and $H$ differs from $G$ in a single edge $uv$, i.e., $H=G-uv$. Then $d_{H}(u)=d_{H}(v)=\Delta-1$. Note that $\Delta\geq k\geq 2$. We consider the following two cases:

{\it Case 1.} $\Delta=2$.

In this case, $G$ must be the cycle $C_{n}$ on $n$ vertices, and thus $H$ is the path $P_{n}$ on $n$ vertices. Further, noticing that $q(P_{n})=2+2\cos\frac{\pi}{n}$ and $\sin x>x-x^{3}/6$, one check that
$$2\Delta-q(H)=2(1-\cos\frac{\pi}{n})=4\sin^{2}\frac{\pi}{2n}> \frac{2}{2n^{2}-7n+9},$$
as desired, completing the proof of Case 1.

{\it Case 2.} $\Delta\geq 3$.

In this case, note that $H$ is connected since $k\geq 2$. Then let $x=(x_{1}, x_{2}, \cdots, x_{n})^{T}$ be the unique unit positive eigenvector of $Q(H)$ corresponding to $q(H)$. Clearly, $x_{1}^{2}+x_{2}^{2}+\cdots+x_{n}^{2}=1$. And let $w$ be a vertex such that $x_{w}=max_{1\leq i\leq n}x_{i}$. By similar arguments as the proof of Theorem \ref{th:t1}, we have that $u\neq w$ and $v\neq w$.
We also find that
\begin{align*}
\displaystyle 2\Delta-q(H)&=2\Delta\sum\limits_{i=1}^{n}x_{i}^{2}-\sum\limits_{i=1}^{n}d_{i}x_{i}^{2}-2\sum\limits_{v_{i}v_{j}\in E(H)}x_{i}x_{j}\\
            &=2(x_{u}^{2}+x_{v}^{2})+\sum\limits_{i=1}^{n}d_{i}x_{i}^{2}-2\sum\limits_{v_{i}v_{j}\in E(H)}x_{i}x_{j}
 \end{align*}
\begin{equation}\label{eq:e6}
~~~~~~=2(x_{u}^{2}+x_{v}^{2})+\sum\limits_{v_{i}v_{j}\in E(H)}(x_{i}-x_{j})^{2}.
 \end{equation}
Since $\kappa(H-v)\geq k-1$, again by Menger's Theorem, there are (at least) $k-1$ vertex-disjoint paths joining $w$ and $u$ in $H-v$, say $P_{1}, P_{2},\cdots,P_{k-1}$, which are as short as possible. Clearly, each of these paths contains only one vertex in $N_{H}(u)$, and then $\sum\limits_{t=1}^{k-1}|V(P_{t})|\leq n-\Delta+3k-5$. Thus by the Cauchy-Schwarz inequality we get
\begin{align*}
\displaystyle \sum\limits_{v_{i}v_{j}\in E(G)}(x_{i}-x_{j})^{2}&\geq \sum\limits_{t=1}^{k-1}\sum\limits_{v_{i}v_{j}\in E(P_{t})}(x_{i}-x_{j})^{2}\\
                    &\geq \sum\limits_{t=1}^{k-1}\frac{1}{|V(P_{t})|-1}(\sum\limits_{v_{i}v_{j}\in E(P_{t})}(x_{i}-x_{j}))^{2}
 \end{align*}

\begin{align*}
                    &=(\sum\limits_{t=1}^{k-1}\frac{1}{|V(P_{t})|-1})(x_{w}-x_{u})^{2}\\
                    &\geq \frac{(k-1)^{2}}{\sum\limits_{t=1}^{k-1}(|V(P_{t})|-1)}(x_{w}-x_{u})^{2}
 \end{align*}
 \begin{equation}\label{eq:e7}
~~~~~~~~~~~~~~~\geq \frac{(k-1)^{2}}{n-\Delta+2k-4}(x_{w}-x_{u})^{2}.
 \end{equation}
Combining (\ref{eq:e6}) and (\ref{eq:e7}), and using Lemma \ref{le:l2}, we have
\begin{align*}
\displaystyle 2\Delta-q(H)&> 2x_{u}^{2}+\frac{(k-1)^{2}}{n-\Delta+2k-4}(x_{w}-x_{u})^{2}
 \end{align*}
 \begin{equation}\label{eq:e8}
~~~~~~~~\geq \frac{2(k-1)^{2}}{2(n-\Delta+2k-4)+(k-1)^{2}}x_{w}^{2}.
 \end{equation}
Let
$$B=\frac{2(k-1)^{2}}{2(n-\Delta)(n-\Delta+2k-4)+(n+1)(k-1)^{2}}.$$
Next we will show that $2\Delta-q(H)>B$.

Suppose that $N_{H}(u)= \{ u_{1}, u_{2}, \cdots, u_{\Delta-1} \}$. Here $w$ may be $u_{t}$ for some $t \in \{1, 2, \cdots, \Delta-1\}$, if this is the case, for convenience, we assume $w=u_{\Delta-1}$.

{\it Subcase 2.1.} $x_{u}^{2}+x_{v}^{2}\geq B/2$.

In this case, from (\ref{eq:e6}), we can get
$$2\Delta-q(H)> 2(x_{u}^{2}+x_{v}^{2})> B.$$

{\it Subcase 2.2.} $\sum\limits_{t=1}^{\Delta-2}x_{u_{t}^{2}}\geq \frac{\Delta}{2}B$

In this case, for avoiding the possible case of $w=u_{\Delta-1}$, then using (\ref{eq:e6}) and Lemma \ref{le:l2}, we obtain
\begin{align*}
\displaystyle 2\Delta-q(H)&\geq 2x_{u}^{2}+\sum\limits_{t=1}^{\Delta-2}(x_{u_{t}}-x_{u})^{2}\\
                   &=\sum\limits_{t=1}^{\Delta-2}[\frac{2}{\Delta-2}x_{u}^{2}+(x_{u_{t}}-x_{u})^{2}]\\
                   &\geq \frac{2}{\Delta}\sum\limits_{t=1}^{\Delta-2}x_{u_{t}}^{2}\geq B.
 \end{align*}

{\it Subcase 2.3.} $x_{u}^{2}+x_{v}^{2}< B/2$ and $\sum\limits_{t=1}^{\Delta-2}x_{u_{t}^{2}}< \frac{\Delta}{2}B$.

In this case, noticing that
$$x_{w}^{2}\geq (1-x_{u}^{2}-x_{v}^{2}-\sum\limits_{t=1}^{\Delta-2}x_{u}(t)^{2})/(n-\Delta)> (1-\frac{\Delta+1}{2} B)/(n-\Delta),$$
and from (\ref{eq:e8}) again, we have
\begin{align*}
\displaystyle 2\Delta-q(H) &> \frac{2(k-1)^{2}}{2(n-\Delta+2k-4)+(k-1)^{2}}x_{w}^{2}\\
             &>\frac{2(k-1)^{2}}{[2(n-\Delta+2k-4)+(k-1)^{2}](n-\Delta)}(1-\frac{\Delta+1}{2} B)\\
             &=\frac{2(k-1)^{2}}{2(n-\Delta)(n-\Delta+2k-4)+(n+1)(k-1)^{2}}=B.
\end{align*}
This completes the proof of Theorem \ref{th:t2}.\quad $\square$

\section{Final remarks}
It is easy to prove that when $k>2\sqrt{\frac{(n-\Delta)(n+\Delta-4)}{n(4D-3)-2}}+1$, bound (\ref{eq:e4}) is always better than bound (\ref{eq:e3}). Indeed, when  $k>2\sqrt{\frac{(n-\Delta)(n+\Delta-4)}{n(4D-3)-2}}+1$, from (\ref{eq:e4}), then we have
\begin{align*}
\displaystyle 2\Delta-q(H)&>\frac{2(k-1)^{2}}{2(n-\Delta)(n+\Delta-4)+(n+1)(k-1)^{2}}\\
&=\frac{1}{(n-\Delta)(n+\Delta-4)/(k-1)^{2}+(n+1)/2}\\
&>\frac{1}{n(D-1/4)}.
\end{align*}
On the other hand, when $k<\frac{2(n-\Delta)}{\sqrt{n(4D-3)-2}}+1$, bound (\ref{eq:e3}) is always better than bound (\ref{eq:e4}). Indeed, when $k<\frac{2(n-\Delta)}{\sqrt{n(4D-3)-2}}+1$, from (\ref{eq:e3}), we have
\begin{align*}
\displaystyle 2\Delta-q(H)&>\frac{1}{n(D-1/4)}\\
&>\frac{1}{(n-\Delta)^{2}/(k-1)^{2}+(n+1)/2}\\
&=\frac{2(k-1)^{2}}{2(n-\Delta)^{2}+(n+1)(k-1)^{2}}\\
&> \frac{2(k-1)^{2}}{2(n-\Delta)(n-\Delta+2k-4)+(n+1)(k-1)^{2}}.
\end{align*}

To provide some preliminary evidence, we here list some values of bounds (\ref{eq:e3}) and (\ref{eq:e4}), as shown in Table \ref{tab:tab1}. Graphs $G_{1}$ and $G_{2}$ are the $3$-regular graphs, as shown in Figure 1.  And $G_{11}$ and $G_{12}$ are the maximal subgraphs of $G_{1}$ and $G_{21}$ is the maximal subgraph of $G_{2}$, respectively, as shown in Figure 2.

\begin{figure}[htb]
\centering
\includegraphics[scale=0.3,width=75mm]{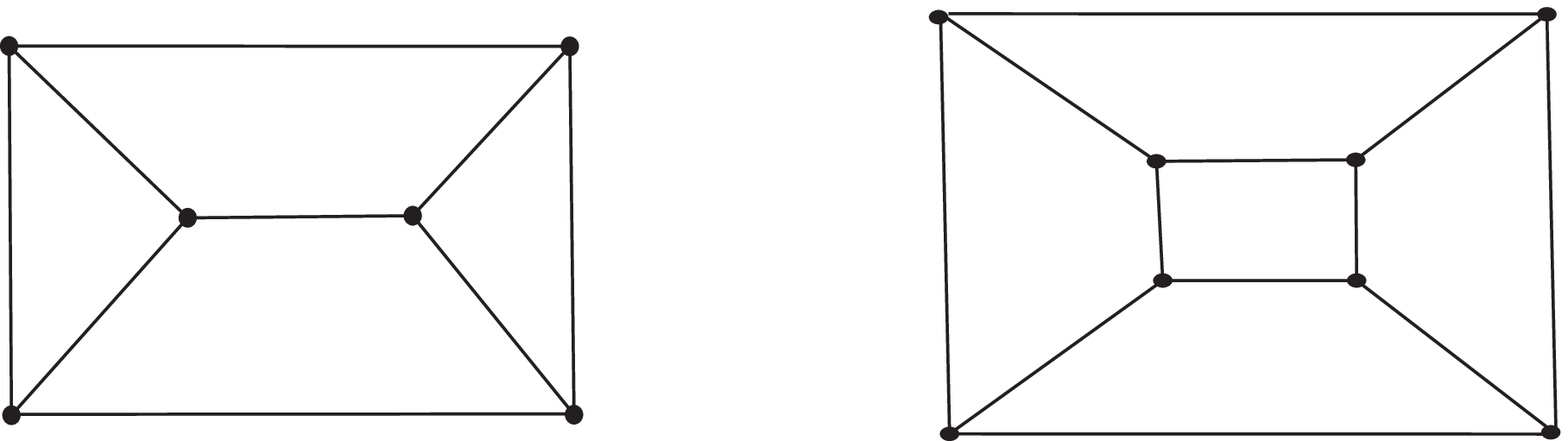}
\caption{3-regular graphs $G_{1}$ and $G_{2}$}
\vskip2mm

\centering
\includegraphics[scale=0.6,width=90mm]{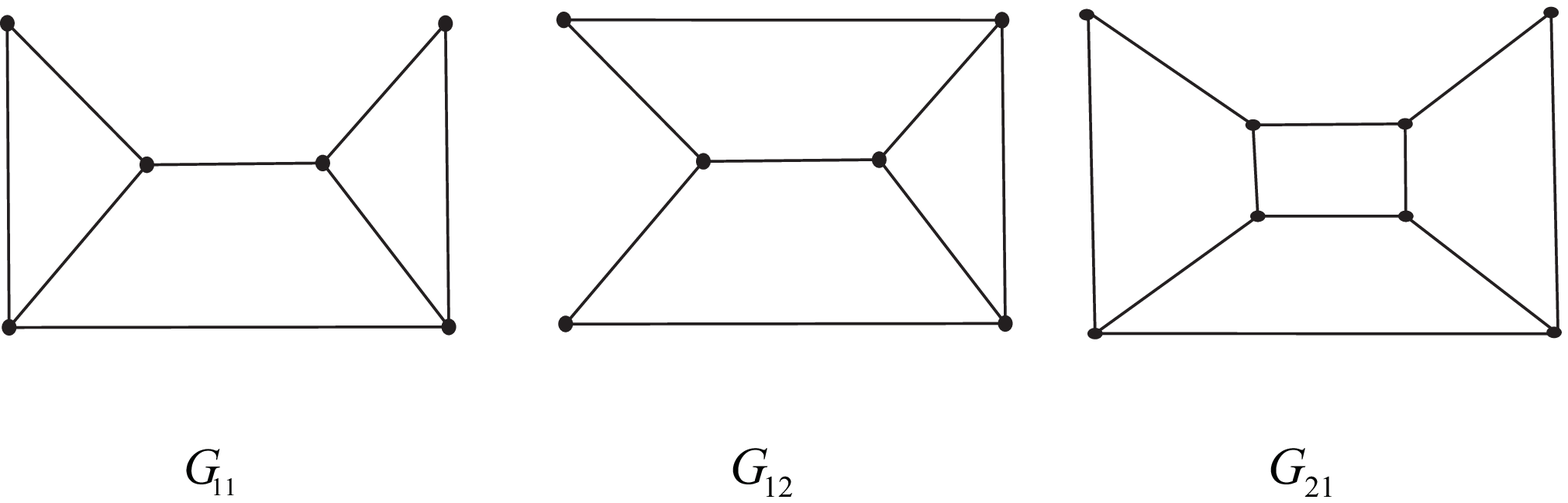}
\caption{The subgraphs of $G_{1}$ and $G_{2}$}
\end{figure}

\begin{table}[htbp]
\centering
\caption{Bounds (\ref{eq:e3}) and (\ref{eq:e4})  of maximal subgraphs $H$ for regular graphs}\label{tab:tab1}
\begin{tabular}{ccccc}
\hline
Graph  & Maximal subgraph $H$  & $2\Delta-q(H)$ & (\ref{eq:e3}) & (\ref{eq:e4}) \\
\hline
$C_{6}$  &  $P_{6}$  &  0.268    &     \textbf{0.0606}      &    0.05128 \\
$C_{12}$  &  $P_{12}$  &  0.0682    &     \textbf{0.0159 }     &    0.0094 \\

$K_{6}$  &  $K_{6}-e$  &  0.5359    &     0.2222      &    \textbf{0.25397} \\

$K_{12}$  &  $K_{12}-e$  &  0.2918    &     0.1111     &    \textbf{0.14948} \\

$G_{1}$  &  $G_{11}$  &  0.4384    &     0.0952      &    \textbf{0.1379} \\

$G_{1}$  &  $G_{12}$  &  0.4113    &     0.0952      &    \textbf{0.2069} \\

$G_{2}$  &  $G_{21}$  &  0.2907    &     0.0714      &    \textbf{0.0816} \\
\hline
\end{tabular}
\end{table}

\end{document}